\documentclass{article}
\usepackage{amssymb}
\usepackage{here}
\usepackage{geometry}
\usepackage{amsfonts}
\usepackage{indentfirst}
\usepackage{mathenv}
\usepackage{amsmath}
\usepackage[english,francais]{babel}
\usepackage[T1]{fontenc}
\usepackage{authblk}
\newtheorem{Theorem}{Theorem}[section]

\newtheorem{Proposition}[Theorem]{Proposition}

\newtheorem{Conjecture}{Conjecture}

\newcommand{\cqfd}
{%
\mbox{}%
\nolinebreak%
\hfill%
\rule{2mm}{2mm}%
\medbreak%
\par%
}
\title{An application of a conjecture due to Ervedoza and Zuazua concerning the observability of the heat equation in small time to a conjecture due to Coron and Guerrero concerning the uniform controllability of a convection-diffusion equation in the vanishing viscosity limit}
\author[1]{Pierre Lissy\footnote{lissy@ann.jussieu.fr}\thanks{Work supported by ERC advanced grant 266907 (CPDENL) of the 7th Research Framework Programme (FP7)}}
\affil[1]{UPMC Univ Paris 06, UMR 7598, Laboratoire Jacques-Louis Lions, F-75005, Paris, France.}
\date{\empty}
\begin{document}
\maketitle
\selectlanguage{english}

\begin{abstract}
The aim of this short paper is to explore a new connection between a conjecture concerning sharp boundary observability estimates for the $1$-D heat equation in small time 
and a conjecture concerning the cost of null-controllability for a $1$-D convection-diffusion equation with constant coefficients controlled on the boundary in the vanishing viscosity limit, in the spirit of what is done in [Pierre Lissy, A link between the cost of fast controls for the $1$-D heat
              equation and the uniform controllability of a 1-D
              transport-diffusion equation, C. R. Math. Acad. Sci. Paris, Volume 352, 2012]. We notably establish that the first conjecture implies the second one as soon as the speed of the transport part is non-negative in the transport-diffusion equation.

							\end{abstract}
	\noindent{\bf Keywords:}
Transport-diffusion equation; Null controllability; Vanishing viscosity limit.
\vspace{2\baselineskip}						
\section{Introduction}
Let us consider some $T>0$, $L>0$, and the following one-dimensional heat equation with boundary control at the left side of the boundary:

\begin{equation}
\left\{
\begin{aligned}
y_t-y_{xx}&=0&\mbox{ in } (0,T)\times (0,L),\\
y(\cdot,0)&=u(t)&\mbox{ in }  (0,T),\\
y(\cdot,L)&=0 &\mbox{ in }  (0,T),
\end{aligned}
\right .
\label{heatez}
\end{equation}
with initial condition $y(0,\cdot)=y^0\in H^{-1}(0,L)$ and control $u\in L^2(0,T)$. It is well-known that this control operator is admissible for initial data in $H^{-1}(0,L)$ and that equation \eqref{heatez} is null-controllable in arbitrary small time, thanks (notably) to the infinite speed of propagation of the information (see \cite{71FR} for the one-dimensional case, or more recently \cite{96FIbook}, \cite{95LR} and \cite{MR2593342} for results in any space dimension with in addition terms of order $0$ and $1$ in the equation). A more challenging (and still widely open, especially in the multi-dimensional case) question is what is usually called \emph{the cost of fast controls}, i.e. what is the energy required to steer the system to $0$ at time $T$ when $T\rightarrow 0$. More precisely, let us call $$\mathcal T_D(y^0,T):=\{(y,u)\in C^0([0,T],H^{-1}(0,L))\times L^2(0,T)|(y,u)\mbox{ verifies \eqref{heatez}, }y(0,\cdot)=y^0\mbox{ and } y(T,\cdot)=0\}$$
and
\begin{equation*}
C_D(T,L)=\sup_{y^0\in H^{-1}(0,L)}\inf_{(y,u)\in\mbox{ }\mathcal T_D(y^0,T)}\frac{||u||_{L^2(0,T)}}{||y^0||_ {H^{-1}(0,L)}}.
\end{equation*}

One can prove (see for example \cite[Section 2.3]{MR2302744}) that $C_D(T,L)$ is always finite. Moreover, for every $y^0\in H^{-1}(0,L)$, the infimum is always reached at a unique $u_{opt}$, which corresponds to the $L^2(0,T)$-projection of the vector $0$ on the subspace of all the controls $u$ that steer $y^0$ to $0$ at time $T$ (i.e. $u_{opt}$ is the control having the smallest norm, which justifies the expression ``cost of the control'' for the quantity $C_D(T,L)$). $C_D(T,L)$ can be seen as the smallest constant $C>0$ such that for every $y^0\in H^{-1}(0,L)$, there exists some control $u$ such that the corresponding solution $y$ of \eqref{heatez} with initial condition $y^0$ verifies $y(T)=0$ and 
\begin{equation}\label{coutH}||u||_{L^2(0,T)}\leqslant C||y^0||_{H^{-1}(0,L)}.
\end{equation}
Let us call 
$$\beta^*:= \underset{T\rightarrow 0}\limsup\mbox{ }T\ln(C_D(T,L))$$
and
$$\beta_*:= \underset{T\rightarrow 0}\liminf\mbox{ }T\ln(C_D(T,L)).$$

It is proved in \cite{MR805273} that $\beta_*>0$ and in \cite{MR743923} that $\beta^*<+\infty$. This means that $C_D(T,L)$ behaves roughly like $e^{\frac{K}{T}}$ for small enough $T$, where $K$ is independent of $T$.
These results have been made more precise later: One has  $\beta_*\geqslant L^2/4$ (see \cite{MR2062431}), which can be notably proved by constructing a singular solution of the heat equation and using a formula of Varadhan concerning the heat kernel in small time (cf. \cite{MR0208191}). The question of upper bounds for $\beta^*$ has also been studied by numerous authors  (see notably \cite{MR743923}, \cite{MR869762}, \cite{2004-Miller} or \cite{MR2370666}) and the best upper bound known is $\beta^*\leqslant 3L^2/4$ as obtained in \cite{TT}. Let us mention here that thanks to the \emph{transmutation} of solutions of wave-type equations into solutions of heat-type equations studied systematically in \cite{MR2246098}, a result for the cost of fast controls in one space dimension for boundary or distributed controls gives a corresponding result in the multi-dimensional case in a (smooth enough) bounded domain $\Omega$ of $\mathbb R^n$ ($n\in\mathbb N^*$) with control domain $\omega$ which can be either a subset of $\partial\Omega$ or distributed into $\Omega$ verifying the geometric control condition GCC of \cite{Bardos-Lebeau-Rauch}, which justifies the specific study of the one-dimensional case. However, it is not natural to impose GCC for the control domain $\omega$ in the case of parabolic systems. Moreover, in this case tone can bound $\beta^*$ from above by a multiple of the length of the longest geodesic of $\overline \Omega$ not intersecting $\omega$, which is not natural for heat-type equations (one expects $\beta^*$ to be some multiple of $\sup_{y\in\Omega} d(y,\omega)$ according for example to \cite[Theorem 2.1]{2004-Miller}).  Using the transmutation method for control domains not verifying GCC is possible (see \cite{MR2833258}), unfortunately it does not seem to provide precise estimates on $C_D(T,L)$.

In fact, it is conjectured (notably in \cite{2004-Miller} and \cite{EZ2011}) that in the one-dimensional case one has exactly
\begin{Conjecture}
\label{Cout-H}
$\beta^*:=(L^2/4)^+$, i.e. for every $L>0$, for every $K>L^2/4$, there exists some $C(K)>0$ such that for every $T$ small enough,
$$C_D(T,L)\leqslant C(K)e^{\frac{K}{T}}.$$
\end{Conjecture}
For corresponding conjectures in the multi-dimensional case, see \cite{2004-Miller} or \cite[Section 5]{EZ2011}.

It is well-known, using the duality between observability and controllability (see for example \cite{77DR}), that the null-controllability of equation \eqref{heatez} with cost $C_D(T,L)$ is exactly equivalent to proving the following inequality, called \emph{observability inequality}:
\begin{equation}\label{ineg-obsez}
\int_0^L |\varphi(T,x)|^2\leqslant C_D(T,L)^2\int_0^T |\partial_x\varphi(t,0)|^2dt,
\end{equation}
for every $\varphi$ satisfying
\begin{equation}
\left\{
\begin{aligned}\label{obsh}
\varphi_t-\varphi_{xx}&=0&\mbox{ in }(0,T)\times(0,L),
\\ \varphi(\cdot,0)&=0&\mbox{ in }(0,T),
\\ \varphi(\cdot,L)&=0&\mbox{ in }(0,T),
\\ \varphi(0,\cdot)&=\varphi^0&\mbox{ in }(0,L),
\end{aligned}
\right . 
\end{equation}
with $\varphi^0\in H^1_0(0,L)$. Inequality \eqref{ineg-obsez} is in general set for the adjoint problem of \eqref{heatez} but one can consider instead the forward problem \eqref{obsh} by changing $t$ into $T-t$.
A classical tool to obtain \eqref{ineg-obsez} is to use \emph{parabolic Carleman estimates} as in \cite{96FIbook} (another possible strategy is to use elliptic Carleman estimates  and the FBI transform as in \cite{95LR}). Such inequalities provide integral estimates in finite or infinite time with singular weight at time $t=0$ (because of the irreversibility of the heat semigroup) on a solution $\varphi$ of \eqref{obsh} under the form
\begin{equation}\label{est-Carl}
\int_0^Te^{-\frac{C(L)}{t}}|\varphi(t,x)|^2dt\leqslant C(T,L)\int_0^T |\partial_x\varphi(t,0)|^2dt
\end{equation}
or 
\begin{equation}\label{est-Carl-inf}
\int_0^\infty e^{-\frac{ C_\infty(L)}{t}}|\varphi(t,x)|^2dt\leqslant C(T,L)\int_0^T |\partial_x\varphi(t,0)|^2dt,
\end{equation}
for some constants $C(L)$, $C_\infty(L)$ (not depending on $T$) and some constants $C(T,L)$. The links between estimates of the form \eqref{ineg-obsez}, \eqref{est-Carl} and \eqref{est-Carl-inf} has been carefully studied in \cite{MR2287707} (it is also worth mentioning here \cite{MR2679651} which concentrates on the link between \eqref{ineg-obsez} and spectral inequalities coming from the Lebeau-Robbiano strategy).

The main problem of Carleman estimates is that they are not really adapted to the context of the cost of fast controls because they are in general not precise enough to ensure that the constant $C$ and $C_\infty$ appearing in such estimates are optimal. It is proved in the very nice paper \cite{EZ2011} (using a kind of reverse transmutation method where the solutions of some wave equation are written in terms of solutions of the corresponding heat equation) that the following estimate in infinite-time horizon holds:

\begin{equation}\label{sharp-est-inf}
\int_0^\infty\int_0^L e^{-\frac{L^2}{2t}}|\varphi(t,x)|^2dxdt\leqslant C(L)\int_0^\infty |\partial_x\varphi(t,0)|^2dt,
\end{equation}
where $C(L)$ might depend on $L$.
Moreover,  this estimate is sharp in the one-dimensional case (as proved in \cite{71FR}) and characterizes the reachable states of equation \eqref{heatez}. The authors are then able to derive from inequality \eqref{sharp-est-inf} the following estimate in finite time:
\begin{equation}\label{sharp-est-f}
\int_0^\infty\int_0^L  e^{-\frac{L^2}{2t}}|\varphi(t,x)|^2dxdt\leqslant C_{int}(T,L)\int_0^T |\partial_x\varphi(t,0)|^2dt,
\end{equation}
where $C_{int}(T,L)$ is a constant that depends on $T$ and $L$. Unfortunately, since \eqref{sharp-est-f} is obtained by using a reasoning by contradiction, the authors where unable to estimate precisely the constant $C_{int}(T,L)$.

A very natural conjecture (cf. \cite[Section 1.2, Section 3.2, Section 5]{EZ2011}) would be that the constant $C_{int}(T,L)$ does not blow up in a too violent way, in the following sense:
\begin{Conjecture}\label{Conj-Int} Let $\delta>0$ and $L>0$. One can choose $C_{int}(T,L)$ such that $$C_{int}(T,L)=\underset{T\rightarrow 0}O(e^{\frac{\delta}{T}}).$$ \end{Conjecture}
Conjecture~\ref{Conj-Int} would notably be true if for example $C_{int}(T,L)$ is some fraction of (some power of) $T$.

Let us mention here the well-known fact that the dissipative character of the parabolic systems implies the following result: 
\begin{Proposition}\label{intimponc}
Conjecture~\ref{Conj-Int} is stronger than Conjecture~\ref{Cout-H}.
\end{Proposition}
For the sake of clarity, we recall here briefly a possible proof, without claim of originality (this is very similar to the proof of the fact that Carleman estimate $\Rightarrow$ observability).

\textbf{Proof of Proposition~\ref{intimponc}.}

 Let us consider some $0<r<1$ (which is destined to be close to $1$). Using the fact that the $L^2$-norm of the solution of \eqref{obsh} is non-increasing, we obtain
$$\int_0^T\int_0^L e^{-\frac{L^2}{2t}}|\varphi(t,x)|^2dxdt\geqslant e^{-\frac{L^2}{2Tr}}\int_{rT}^{T}\int_0^L |\varphi(t,x)|^2dxdt\geqslant (1-r)Te^{-\frac{L^2}{2rT}}\int_0^L |\varphi(T,x)|^2dx.$$
For a given $K>L^2/4$, thanks to Conjecture~\ref{Conj-Int}, there exists some constant $C>0$ and some $r<1$ (that might depend on $L$ but not on $T$) such that for every $T>0$ small enough one has 
$$\frac{e^{\frac{L^2}{2rT}}}{T(1-r)}C_{int}(T,L)\leqslant Ce^{\frac{2K}{T}}.$$  We deduce that  
$$\int_0^L |\varphi(T,x)|^2dx\leqslant Ce^{\frac{2K}{T}}\int_0^T |\partial_x\varphi(t,0)|^2dt,$$
which gives exactly Conjecture~\ref{Cout-H}.
\cqfd

The goal of what follows is to explain how Conjecture~\ref{Conj-Int} can be linked to another famous conjecture stated in \cite{05aa} concerning the uniform controllability of a transport-diffusion equation in one space dimension with constant coefficients in the vanishing viscosity limit, in the spirit of what was done by the author in \cite{MR2956149}.

Let us consider some constant $M\not =0$ (supposed to be independent of $x$ and $t$) and some viscosity coefficient $\varepsilon\in (0,1)$ (which is destined to tend to $0$). We are interested in the following family of transport-diffusion equations
\begin{equation}
\left\{
\begin{aligned}\label{tdez}
y_t-\varepsilon y_{xx}+My_x&=0&\mbox{ in } (0,T)\times (0,L),&\\
y(\cdot,0)&=v(t)&\mbox{ in }  (0,T),&\\
y(\cdot,L)&=0 &\mbox{ in }  (0,T),&
\end{aligned}
\right . 
\end{equation}
with initial condition $y^0\in H^{-1}(0,L)$ and control $v\in L^2(0,T)$ and $\varepsilon \in (0,1)$.
If $\varepsilon$ is taken equal to $0$ and if the initial condition $y^0$ is taken in $L^2(0,T)$, we obtain a transport equation at constant speed $M$
\begin{equation}
\left\{
\begin{aligned}\label{transp}
y_t+My_x&=0&\mbox{ in } (0,T)\times (0,L),&\\
y(\cdot,0)&=v(t)&\mbox{ in }  (0,T),&\\
y(\cdot,L)&=0 &\mbox{ in }  (0,T),&
\end{aligned}
\right . 
\end{equation}
which is known to be null-controllable if and only if $T\geqslant L/|M|$, the optimal control in $L^2$-norm is in this case the null function (see for example \cite[Section 2.1]{MR2302744}).
Let us define 
$$\mathcal T_{TD}(y^0,T):=\{(y,u)\in C^0([0,T],H^{-1}(0,L))\times L^2(0,T)|(y,u)\mbox{ verifies \eqref{tdez}, }y(0,\cdot)=y^0\mbox{ and } y(T,\cdot)=0\}$$
and
\begin{equation*}
C_{TD}(T,L)=\sup_{y^0\in H^{-1}(0,L)}\inf_{(y,u)\in\mbox{ }\mathcal T_{TD}(y^0,T)}\frac{||u||_{L^2(0,T)}}{||y^0||_ {H^{-1}(0,L)}}.
\end{equation*}

Since one can prove (see \cite[Appendix A]{05aa}) that the solution of \eqref{tdez} with initial condition $y^0\in L^2(0,L)$ converges in some sense to the one of \eqref{transp} when $\varepsilon\rightarrow 0$, one might reasonably expect that $C_{TD}(T,L,M,\varepsilon)\rightarrow +\infty$ for  $T<L/|M|$ and $C_{TD}(T,L,M,\varepsilon)\rightarrow 0$ for $T>L/|M|$.  

However, it is proved in \cite{05aa}, as expected, that one has $$C_{TD}(T,L,M,\varepsilon)\geqslant Ce^{\frac{K}{\varepsilon}}$$ for some constants $C,K$ independent of $\varepsilon$ if $T< L/M$ for $M>0$, but what is unexpected is that $$C_{TD}(T,L,M,\varepsilon)\geqslant Ce^{\frac{K}{\varepsilon}}$$ for some $C,K$ independent of $\varepsilon$ if $T< 2L/|M|$ for $M<0$. This surprising result led the authors to make the following conjecture concerning positive results for the uniform controllability of the family of equations \eqref{tdez} in large time, which is still not decided to be true or false:
\begin{Conjecture}Let $T>0$, $L>0$ and $M\not =0$ be given. Then $C_{TD}(T,L,M,\varepsilon)\rightarrow 0$ as $\varepsilon\rightarrow 0^+$ as soon as $T>L/M$ for $M>0$ and $T>2L/|M|$ for $M<0$. \label{conjCGez}
\end{Conjecture}
In \cite{05aa}, it is proved the exponential decay of the cost of the control when $\varepsilon\rightarrow 0^+$ for sufficiently large times, namely $T>4.3 L/M$ (resp.  $T>57.2L/|M|$) if $M>0$ (resp. $M<0$), which was extended to varying in time and space (and regular enough) speed $M$  and arbitrary space dimension in \cite{2005-Guerrero-Lebeau}. In both articles \cite{05aa} and \cite{2005-Guerrero-Lebeau}, the strategy is to derive a Carleman estimate which takes into account the transport term and then use dissipation result adapted to the equation similar to what is done in \cite{MR1472116}. In \cite{2005-Guerrero-Lebeau}, the authors also obtained an exponential growth of the cost of the control as $\varepsilon\rightarrow 0^+$ for small times. 

The upper bounds concerning  the uniform controllability of equation \eqref{tdez} 
have been improved in \cite{2010-Glass-JFA} which proved the uniform controllability for $T>4.2L/M$ (resp.  $T>6.1L/|M|$) if $M>0$ (resp. $M<0$) by using a method similar to the moment method on the adjoint system of \eqref{tdez} with a well-chosen complex multiplier. The latest improvement were done by the author in \cite{MR2956149}, where the uniform controllability is proved for $T>2\sqrt{3}L/M>3.45 L/M$ (resp. $T>(2\sqrt{3}+2)>5.45L/|M|$) if $M>0$ (resp. $M<0$), which is still quite far from Conjecture~\ref{conjCGez}. However, it is also given a strategy to improve the critical times up to $T>2$ (resp. $T>4$) if $M>0$ (resp. $M<0$) by finding a link between this problem of uniform controllability and Conjecture~\ref{Cout-H} for the heat equation \eqref{heatez} in small time.

We now state our result.
\begin{Theorem}\label{cras-am}
Let $T>0$, $L>0$ and $M>0$ be fixed. Assume that Conjecture~\ref{Conj-Int} is verified. Then there exists some constants $C,K>0$ (independent of $\varepsilon$ but depending possibly on $T>0$, $L>0$ and $M>0$) such that for every $\varepsilon\in (0,1)$,
$$C_{TD}(T,L,\varepsilon,M)\leqslant Ce^{-\frac{K}{\varepsilon}}$$
as soon as
\begin{enumerate}
\item $T> L/M$ for $M>0$,
\item $T> (1+\sqrt{2})L/|M|$ for $M<0$.
\end{enumerate}
\end{Theorem}

Theorem~\ref{cras-am} is quite surprising, because it seems to indicate that Conjecture~\ref{Conj-Int} is much stronger than Conjecture~\ref{Cout-H}, which is unexpected because of one may think that what happens near $t=+\infty$ is quite negligible compared to what happens near $t=T$. Moreover, we think that Theorem~\ref{cras-am} is of interest because it will enable people to solve both Conjecture~\ref{Cout-H} and Conjecture~\ref{conjCGez} in the case $M>0$ simultaneously. Theorem~\ref{cras-am} suggests that it is worth trying to concentrate on proving Conjecture~\ref{conjCGez} instead of proving Conjecture~\ref{Cout-H}.

\section{Proof of Theorem~\ref{cras-am}}

Let $\psi^0\in H^1_0(0,L)$ and let $\psi$ be the solution of the following forward problem:
\begin{equation}
\left\{
\begin{aligned}\label{tdezret}
\psi_t-\varepsilon \psi_{xx}-M\psi_x&=0&\mbox{ in } (0,T)\times (0,L),&\\
\psi(\cdot,0)&=0&\mbox{ in }  (0,T),&\\
\psi(\cdot,L)&=0&\mbox{ in }  (0,T),&\\
\psi(0,\cdot)&=\psi^0&\mbox{ in }  (0,T).&
\end{aligned}
\right . 
\end{equation}

We use the same kind of transformation than in \cite[Proof of Lemma 2.1]{MR2956149}, and we call
\begin{gather}\label{chvar}\varphi(t,x):=e^{\frac{M^2t}{4\varepsilon^2}+\frac{Mx}{2\varepsilon}}\psi(\frac{t}{\varepsilon},x).\end{gather}
Then 
\begin{equation*}\varphi_t(t,x)- \varphi_{xx}(t,x)=e^{\frac{M^2t}{4\varepsilon^2}+\frac{Mx}{2\varepsilon}}(\frac{M^2}{4\varepsilon^2}\psi(\frac{t}{\varepsilon},x)+\frac{\psi(\frac{t}{\varepsilon},x)}{\varepsilon}-\frac{M^2}{4\varepsilon^2}\psi(\frac{t}{\varepsilon},x)-\psi(\frac{t}{\varepsilon},x)-\frac{M}{\varepsilon}\psi_x)=0.\end{equation*}
Hence $\varphi$ is a solution of \eqref{obsh} on the interval $(0,\varepsilon T)\times (0,L)$, with initial data $\varphi^0(x)=e^{\frac{Mx}{2\varepsilon}}\psi^0(x)$. 
Since $\varepsilon\rightarrow 0$, this means that we now work in small time and the dependance in $\varepsilon$ now only appears in the time variable.
Moreover, one has 
 \begin{gather}\label{dchvar}\partial_x \varphi(t,0)=\frac{M}{2\varepsilon}e^{\frac{M^2t}{4\varepsilon^2}}\psi(\frac{t}{\varepsilon},0)+e^{\frac{M^2t}{4\varepsilon^2}}\partial_x \psi(\frac{t}{\varepsilon},0)=e^{\frac{M^2t}{4\varepsilon^2}}\partial_x \psi(\frac{t}{\varepsilon},0).\end{gather}

Let us consider some $0<a<1$ (which is destined to be close to $1$). The map $t\mapsto e^{-\frac{L^2}{2t}}$ is increasing on $(0,\infty)$, so we deduce that 
\begin{equation}\label{coupea}
\left .
\begin{aligned}
e^{-\frac{L^2}{2a\varepsilon T}}\int_{\varepsilon aT}^{\varepsilon T}\int_0^L|\varphi(t,x)|^2dxdt&\leqslant \int_{\varepsilon aT}^{\varepsilon T}\int_0^Le^{-\frac{L^2}{2t}}|\varphi(t,x)|^2dxdt
\\ \mbox{ }&\leqslant \int_{0}^{+\infty}\int_0^Le^{-\frac{L^2}{2t}}|\varphi(t,x)|^2dxdt.
\end{aligned}
\right .
\end{equation}
Since $\varphi$ is a solution of \eqref{obsh}, $t\geqslant 0\mapsto ||\varphi(t,.)||^2_{L^2(0,L)}$ is decreasing so 
\begin{gather}\label{decrphi} (1-a)\varepsilon T \int_0^L |\varphi(\varepsilon T,x)|^2dx\leqslant \int_{\varepsilon aT}^{\varepsilon T}\int_0^L|\varphi(t,x)|^2dxdt.\end{gather}
Using \eqref{coupea} and \eqref{decrphi}, we obtain
\begin{gather}\label{inttoponc} (1-a)\varepsilon T e^{-\frac{L^2}{2a\varepsilon T}}\int_0^L |\varphi(\varepsilon T,x)|^2dx\leqslant \int_{0}^{+\infty}\int_0^Le^{-\frac{L^2}{2t}}|\varphi(t,x)|^2dxdt.\end{gather}
Applying inequality \eqref{sharp-est-f} at time $b\varepsilon T$  with some $b>0$ (which is destined to be close to $0$) and inequality \eqref{inttoponc}, we deduce 
\begin{gather}\label{estphi}(1-a)T e^{-\frac{L^2}{2a\varepsilon T}}\int_0^L |\varphi(\varepsilon T,x)|^2dx \leqslant C_{int}(\varepsilon bT)\varepsilon e^{\frac{bM^2T}{2\varepsilon}}\int_0^{b\varepsilon T} |\partial_x \varphi(t,0)|^2dt. \end{gather}
We now use equalities \eqref{chvar} and \eqref{dchvar} together with inequality \eqref{estphi} to obtain

\begin{equation}\label{intez1}e^{-\frac{L^2}{2a\varepsilon T}+\frac{M^2 T}{2\varepsilon}}\int_0^Le^{\frac{Mx}{\varepsilon}}|\psi(T,x)|^2dx\leqslant C_{int}(\varepsilon bT,L)\int_0^{\varepsilon bT} e^{\frac{M^2t}{2\varepsilon^2}}|\partial_x \psi(\frac{t}{\varepsilon},0)|^2dt.\end{equation}
From \eqref{intez1} we deduce by changing $t$ into $\varepsilon t$ in the integral appearing in the right-hand side
\begin{gather}\label{inegint} e^{-\frac{L^2}{2a\varepsilon T}+\frac{M^2a T}{2\varepsilon}}\varepsilon(1-a)T\int_0^L |e^{\frac{Mx}{2\varepsilon}}\psi(T,x)|^2dxdt\leqslant C_{int}(\varepsilon bT)\varepsilon e^{\frac{bM^2T}{2\varepsilon}}\int_0^{bT} |\partial_x \psi(t,0)|^2dt.\end{gather}
Let us treat separately the cases $M>0$ and $M<0$.
\begin{enumerate}
\item Assume that $M>0$. In this case, inequality \eqref{inegint} provides 
\begin{equation}\label{estez2}\int_0^L |\psi(T,x)|^2dxdt\leqslant \frac{e^{\frac{L^2}{2a\varepsilon T}-\frac{M^2T}{2\varepsilon}+\frac{bM^2T}{2\varepsilon}}}{(1-a)T}C_{int}(\varepsilon b T,L) \int_0^{T} |\partial_x \psi(t,0)|^2dt.\end{equation}
Assume now that $T>0$, $L>0$ and $M>0$ are such that $T>L/M$. Then one obtains 
\begin{equation}\label{unif1}\frac{L^2}{2aT}-\frac{M^2T}{2}+\frac{bM^2T}{2}<0\end{equation}
as soon as $b$ is chosen close enough to $0$ and $a$ is chosen close enough to $1$: The roots of the polynomial (in the variable $T$) 
$$\frac{L^2}{a}+(b-1)M^2 T^2$$ 
converge when $a\rightarrow 1^-$ and $b\rightarrow 0^+$ to the roots of the polynomial 
$$L^2-M^2 T^2,$$ 
which are precisely $-L/M$ and $L/M$. Thanks to \eqref{estez2}, \eqref{unif1} and Conjecture~\ref{Conj-Int}, one deduces that there exists some constants $C,K>0$ independent of $\varepsilon$ (but depending possibly on $L$, $T$ and $M$) such that  
$$C_{TD}(T,L,\varepsilon,M)\leqslant Ce^{-\frac{K}{\varepsilon}}.$$
The first part of Theorem~\ref{cras-am} is proved.
\item Assume now that  $M<0$. In this case, \eqref{inegint} becomes 
\begin{equation}\label{estez3}\int_0^L |\psi(T,x)|^2dxdt\leqslant \frac{e^{\frac{L^2}{2a\varepsilon T}-\frac{M^2 T}{2\varepsilon}+\frac{bM^2T}{2\varepsilon}+\frac{|M|L}{\varepsilon}}}{(1-a)T}C_{int}(\varepsilon b T,L)\int_0^{T} |\partial_x \psi(t,0)|^2dt.\end{equation}
Assume now that $T>0$, $L>0$ and $M<0$ are such that $T>(1+\sqrt{2})L/|M|$. Then one obtains 
\begin{equation}\label{unif2}\frac{L^2}{2aT}-\frac{M^2 T}{2}+\frac{bM^2T}{2}+|M|L<0\end{equation}
as soon as $b$ is chosen close enough to $0$ and $a$ is chosen close enough to $1$: The roots of the polynomial (in the variable $T$) 
$$\frac{L^2}{a}+(b-1)M^2 T^2+2b|M|LT$$ 
converge when $a\rightarrow 1^-$ and $b\rightarrow 0^+$ to the roots of the polynomial 
$$L^2-M^2 T^2+2|M|LT,$$ 
which are precisely $(1-\sqrt{2})L/|M|$ and $(1+\sqrt{2})L/|M|$.

Thanks to \eqref{estez3}, \eqref{unif2} and Conjecture~\ref{Conj-Int}, one deduces that there exists some constants $C,K>0$ independent of $\varepsilon$ (but depending possibly on $L$, $T$ and $M$) such that  
$$C_{TD}(T,L,\varepsilon,M)\leqslant Ce^{-\frac{K}{\varepsilon}},$$ which ends the proof of Theorem~\ref{cras-am}.
\end{enumerate}
\cqfd
\section{Further remarks}
\begin{enumerate}
\item Since one can recover Conjecture~\ref{conjCGez} from Conjecture~\ref{Conj-Int} for $M>0$, maybe $T>2L/|M|$ is not the right conjecture for the minimal time ensuring the uniform controllability for $M<0$ and might be replaced by what is found in Theorem~\ref{cras-am}, i.e. $T>(1+\sqrt{2})L/|M|$.
\item The fact that the integral is on $(0,\infty)$ in the left-hand side of \eqref{sharp-est-f} is crucial in our proof, because it enables us to consider larger times in the left-side ($a$ close to $1$) than in the right-hand side ($b$ close to $0$).
Another possible Conjecture (weaker than Conjecture~\ref{Conj-Int}) would be the following. 
 Let $\delta>0$ and $L>0$. One can find $C_{fin}(T,L)$ such that $$C_{fin}(T,L)=\underset{T\rightarrow 0}O(e^{\frac{\delta}{T}})$$ and
\begin{equation}\label{other-conj}
\int_0^T\int_0^L e^{-\frac{L^2}{2t}}|\varphi(t,x)|^2dxdt\leqslant C_{fin}(T,L)\int_0^T |\partial_x\varphi(t,0)|^2dt.
\end{equation}
However, one can verify that we exactly recover the results of \cite{MR2956149} and notably we do not find constants as good as in Theorem~\ref{cras-am}. This is consistent with the fact that \eqref{other-conj} is roughly equivalent to \eqref{ineg-obsez}, as explained in \cite{MR2287707}.
\end{enumerate}
\bibliographystyle{plain}
\bibliography{Biblio}
\end{document}